\documentclass[amstex,12pt, amssymb]{article}

\usepackage{mathtext}
\usepackage[cp1251]{inputenc}
\usepackage[T2A]{fontenc}
\usepackage[dvips]{graphicx}
\usepackage{amsmath}
\usepackage{amssymb}
\usepackage{amsxtra}
\usepackage{latexsym}
\usepackage{ifthen}

\newcounter{lemma}[section]

\newcounter{corollary}[section]

\newcounter{remark}[section]

\newcounter{theorem}[section]

\newcounter{proposition}[section]

\newcounter{example}

\numberwithin{equation}{section}

\begin{document}

\markboth{V. DESYATKA, E.~SEVOST'YANOV}{\centerline{ON ISOLATED
SINGULARITIES ...}}

\def\cc{\setcounter{equation}{0}
\setcounter{figure}{0}\setcounter{table}{0}}

\overfullrule=0pt


\author{VICTORIA DESYATKA, EVGENY SEVOST'YANOV}

\title{
{\bf ON SINGULARITIES OF MAPPINGS WITH A FINITE LENGTH DISTORTION}}

\date{\today}
\maketitle

\begin{abstract}
We study the possibility of a continuous extension of a class of
mappings to an isolated point on the boundary of a domain. We show
that if some characteristic of this mapping is integrable on almost
all spheres in the neighborhood of at least one point of the
corresponding cluster set, then this mapping has a continuous
extension to the specified point. In particular, this assertion is
true if the specified characteristic is simply Lebesgue integrable
in the neighborhood of at least one limit point.
\end{abstract}


\section{Introduction}

The paper is devoted to the problem of removability of an isolated
singularity in a certain class of mappings. In particular, we are
talking about mappings that satisfy distortion estimates with
respect to the modulus of families of paths; see, for example,
\cite{Ad}, \cite{Af}, \cite{Cr$_1$}--\cite{Cr$_3$}, \cite{MRSY$_2$},
\cite{PSS} and \cite{SalSt}. Note that for essentially singular
points of a quasiregular mapping, the Sokhotsky-Casorati-Weierstrass
theorem is valid, which may be formulated in the following form
(see, e.g., \cite[Corollary~2.11.III]{Ri}).

\medskip
{\bf Theorem~A.} {\it\, Let $b$ be an isolated essential singularity
of a quasimeromorphic mapping $f:G\rightarrow \overline{{\Bbb
R}^n}.$ Then there exists an ${\mathcal F}_{\sigma}$ set $E\subset
\overline{{\Bbb R}^n}$ such that $N(y, f, U\setminus\{b\})=\infty$
for every $y\in \overline{{\Bbb R}^n}\setminus E$ and all such
neighborhoods $U$ of $b.$}

\medskip
Here and below
\begin{equation}\label{eq23}
N(y, f, D)\,=\,{\rm card}\,\left\{x\in D: f(x)=y\right\}\,.
\end{equation}
Let us now look at Theorem~A and the problem of removability of
singularities from the following point of view. It is known that
quasiregular mappings $f:D\setminus\{x_0\}\rightarrow {\Bbb R}^n,$
$x_0\in D,$ $n\geqslant 2,$ satisfy the condition
\begin{equation}\label{eq24}
M(\Gamma)\leqslant \int\limits_{f(D)} K_O\cdot N(y, f,
D\setminus\{x_0\})\cdot \rho_*^n(y)\, dm(y)
\end{equation}
for any function $\rho_*\in{\rm adm}\, f(\Gamma)$ and for any family
$\Gamma$ of paths $\gamma$ in a domain $D,$ where $M$ is a conformal
modulus of families of paths, $K={\rm ess \sup}\, K_O(x, f),$
$K_O(x, f)=K_{O, p}(x, f),$
$$K_{O, p}(x,f)\quad =\quad \left\{
\begin{array}{rr}
\frac{\Vert f^{\,\prime}(x)\Vert^p}{|J(x,f)|}, & J(x,f)\ne 0,\\
1, & f^{\,\prime}(x)=0, \\
\infty, & {\rm in\,other\,cases}
\end{array}
\right.\,,$$
$$
\Vert f^{\,\prime}(x)\Vert\,=\,\max\limits_{h\in {\Bbb R}^n
\backslash \{0\}} \frac {|f^{\,\prime}(x)h|}{|h|}\,,\quad
J(x,f)=\det f^{\,\prime}(x)\,,$$ see \cite[Remarks~2.5.II]{Ri}. Then
Theorem~A implies that, at a minimum, the function $Q(y)=K_O\cdot
N(y, f, D\setminus\{x_0\})$ cannot be integrable in any neighborhood
$U$ of any point $y_0\in {\Bbb R}^n.$ Let us now pose the problem as
follows. Let some open discrete mapping
$f:D\setminus\{x_0\}\rightarrow {\Bbb R}^n,$ $n\geqslant 2,$ $x_0\in
D,$ instead of relation~(\ref{eq24}), satisfies a condition of the
form
\begin{equation}\label{eq24}
M_p(\Gamma)\leqslant \int\limits_{f(D)} Q(y)\cdot \rho_*^p(y)\,
dm(y)\,, \quad n\geqslant p<\infty\,,
\end{equation}
for any function $\rho_*\in{\rm adm}\, f(\Gamma)$ and for any family
$\Gamma$ of paths $\gamma$ in a domain $D.$ {\bf If $x_0$ is an
essential singular point of the mapping $f,$ then what can be said
about the integrability of the function $Q(y)$ in the neighborhood
$U$ of an arbitrary point $y_0\in {\Bbb R}^n$?}

This problem will be further answered for some classes of mappings
satisfying condition~(\ref{eq24}). In this paper, we focus on
mappings with finite length distortion, which are currently one of
the broadest classes of mappings, including quasiconformal,
quasiregular mappings, and some other well-known classes. In this
regard, we should cite the definition of these mappings proposed by
O.~Martio in collaboration with V.~Ryazanov, U.~Srebro, and
E.~Yakubov; see, e.g., \cite{MRSY$_1$}, \cite[Ch.~8]{MRSY$_2$}.

\medskip
Recall some definitions. A Borel function $\rho:{\Bbb
R}^n\,\rightarrow [0,\infty] $ is called {\it admissible} for the
family $\Gamma$ of paths $\gamma$ in ${\Bbb R}^n,$ if the relation
\begin{equation}\label{eq1.4}
\int\limits_{\gamma}\rho (x)\, |dx|\geqslant 1
\end{equation}
holds for all (locally rectifiable) paths $ \gamma \in \Gamma.$ In
this case, we write: $\rho \in {\rm adm} \,\Gamma .$ Let $p\geqslant
1,$ then {\it $p$-modulus} of $\Gamma $ is defined by the equality
\begin{equation}\label{eq1.3gl0}
M_p(\Gamma)=\inf\limits_{\rho \in \,{\rm adm}\,\Gamma}
\int\limits_{{\Bbb R}^n} \rho^p (x)\,dm(x)\,.
\end{equation}
We set $M(\Gamma):=M_n(\Gamma).$

\medskip
We say that a property $P$ holds for {\it $p$-almost every
($p$-a.e.)} paths $\gamma$ in a family $\Gamma$ if the subfamily of
all paths in $\Gamma $, for which $P$ fails, has $p$--mo\-du\-lus
zero.

If $\gamma :\Delta\rightarrow{\Bbb R}^n$ is a locally rectifiable
path, then there is the unique nondecreasing length function
$l_{\gamma}$ of $\Delta$ onto a length interval $\Delta
_{\gamma}\subset{\Bbb R}$ with a prescribed normalization $l
_{\gamma}(t_0)=0\in\Delta _{\gamma},$ $t_0\in\Delta,$ such that $l
_{\gamma}(t)$ is equal to the length of the subpath $\gamma
|_{[t_0,t]}$ of $\gamma$ if $t>t_0,$ $t\in\Delta ,$ and $l
_{\gamma}(t)$ is equal to minus length of $\gamma |_{[t,t_0]}$ if
$t<t_0,$ $t\in\Delta .$ Let $g: |\gamma |\rightarrow{\Bbb R}^n$ be a
continuous mapping, and suppose that the path $\widetilde{\gamma}
=g\circ\gamma$ is also locally rectifiable. Then there is a unique
non--decreasing function $L_{\gamma ,g}: \Delta
_{\gamma}\rightarrow\Delta _{\widetilde{\gamma}}$ such that
$L_{\gamma ,g}\left(l_{\gamma}(t)\right) =
l_{\widetilde{\gamma}}(t)$ for all $t\in\Delta.$ A path $\gamma$ in
$D$ is called here a (whole) {\it lifting} of a path
$\widetilde{\gamma}$ in ${\Bbb R}^n$ under $f:D\rightarrow {\Bbb
R}^n$ if $\widetilde{\gamma} = f\circ\gamma.$

\medskip
We say that a mapping $f:D\rightarrow{\Bbb R}^n$ satisfies the {\it
$L$-pro\-pe\-r\-ty with respect to $(p, q)$-mo\-du\-lus}, iff the
following two conditions hold:

$\left(L^{(1)}_p\right):$ for $p$-a.e. path $\gamma$ in $D,$
$\widetilde{\gamma}=f\circ\gamma$ is locally rectifiable and the
function $L_{\gamma ,f}$ has the $N$-pro\-per\-ty;

$\left(L^{(2)}_q\right):$ for $q$-a.e. path $\widetilde{\gamma}$ in
$f(D),$ each lifting $\gamma$ of $\widetilde{\gamma}$ is locally
rectifiable and the function $L_{\gamma ,f}$ has the
$N^{-1}$-pro\-per\-ty.

The notation $(p,q)$-modulus mentioned above means that the
$L^{(1)}_p$ and $L^{(2)}_q$ properties hold with respect to $p$ and
$q$-moduli, respectively.

\medskip
A mapping $f:D\rightarrow{\Bbb R}^n$ is called a mapping with {\it
finite length $(p, q)$-dis\-tor\-tion}, if $f$ is differentiable
a.e. in $D,$ has $N$-- and $N^{\,-1}$--pro\-per\-ties, and
$L$--pro\-pe\-r\-ty with respect to $(p, q)$-- mo\-du\-lus. The
mappings of finite length $(p, q)$-dis\-tor\-tion are natural
generalization of the mappings with finite length distortion
introduced in the paper \cite{MRSY$_1$}, see also monograph
\cite{MRSY$_2$}.

\medskip
\begin{remark}\label{rem2} Observe that, a mapping $f:D\rightarrow {\Bbb R}^n$ with finite
length $(p, q)$-dis\-tor\-ti\-on satisfies the inequality
\begin{equation}\label{eq5.4AA}
M_p(\Gamma)\leqslant \int\limits_{f(E)}K_{I,p}\left(y, f^{\,-1},
E\right)\cdot\rho_*^p(y)dm(y) \end{equation} for every measurable
set $E\subset D,$ every family $\Gamma$ of paths  $\gamma$ in $E$
and every function $\rho_*(y)\in {\rm adm\,}f(\Gamma),$ where
\begin{equation}\label{eq1U}
K_{I, p}\left(y, f^{\,-1}, E\right):=\sum\limits_{x\in E\cap
f^{\,-1}(y)}K_{O, p} (x,f)\,,\end{equation}
see e.g. \cite[Theorem~1.1]{SalSev}, cf.
\cite[Theorem~8.5]{MRSY$_2$}. Observe that, in this case, $f$ also
satisfies the inequality
\begin{equation}\label{eq2A*}
M_p(\Gamma_f(y_0, r_1, r_2))\leqslant \int\limits_{f(D)\cap
A(y_0,r_1,r_2)} Q(y)\cdot \eta^p (|y-y_0|)\, dm(y)
\end{equation}
with $Q(y)=K_{I, p}\left(y, f^{\,-1}, D\right)$ for $y\in f(D),$
$Q(y)\equiv 0$ for $y\not\in f(D),$ for any Lebesgue measurable
function $\eta: (r_1,r_2)\rightarrow [0,\infty ]$ such that
$\int\limits_{r_1}^{r_2}\eta(r)\, dr\geqslant 1\,.$ Indeed,
let~(\ref{eq2*A}) holds. Now, we put $\rho_*(y):=\eta(|y-y_0|)$ for
$y\in A(y_0,r_1,r_2)\cap f(D),$ and $\rho_*(y)=0$ otherwise. By the
Luzin theorem, we may assume that the function $\rho_*$ is Borel
measurable (see, e.g., \cite[Section~2.3.6]{Fe}). Then,
by~\cite[Theorem~5.7]{Va} we have that
$$\int\limits_{\gamma_*}\rho_*(y)\,|dy|\geqslant
\int\limits_{r_1}^{r_2}\eta(r)\,dr\geqslant 1$$
for any (locally rectifiable) path $\gamma_*\in \Gamma(S(y_0, r_1),
S(y_0, r_2), A(y_0, r_1, r_2)).$ Substituting the function $\rho_*$
in~(\ref{eq5.4AA}), we obtain the desired ratio~(\ref{eq2A*}).
\end{remark}

\medskip
We set
\begin{equation}\label{eq5} C(x, f):=\{y\in \overline{{\Bbb
R}^n}:\exists\,x_k\in D: x_k\rightarrow x, f(x_k) \rightarrow y,
k\rightarrow\infty\}\,.
\end{equation}

\medskip
The following result holds.

\medskip
\begin{theorem}\label{th2}
{\it\, Let $D$ be a domain in ${\Bbb R}^n,$ $n\geqslant 2,$ let
$p\geqslant n$ and let $x_0\in D.$ Let
$f:D\setminus\{x_0\}\rightarrow {\Bbb R}^n$ be an open discrete
mapping with finite length $(p, q)$-dis\-tor\-ti\-on. Let
$Q(y):=K_{I, p}\left(y, f^{\,-1}, D\setminus\{x_0\}\right)\geqslant
c$ for some $c>0$ and for a.e $y\in {\Bbb R}^n.$ If $x_0$ is an
essential singularity of $f,$ then $K_{I, p}\left(y, f^{\,-1},
D\setminus\{x_0\}\right)\not\in L^1(U)$ for any neighborhood $U$ of
$y_0$ and for any $y_0\in C(x_0, f).$

\medskip
In particular, if $K_{I, p}\left(y, f^{\,-1}, D\right)\leqslant
Q(y)\cdot N(y, f, D\setminus\{x_0\})$ for a.e. $y\in
f(D\setminus\{x_0\})$ and some Lebesgue measurable function $Q:{\Bbb
R}^n\rightarrow [0, \infty]$ such that $Q\in L^1(U),$ then $N(y, f,
V\setminus\{x_0\})$ is unbounded for any neighborhood $V$ of $x_0.$
}
\end{theorem}

\medskip
\begin{remark}\label{rem1A}
We should clarify that Theorem~\ref{th2} is a version of Theorem~A
for mappings with finite length $(p, q)$-distortion: in both cases,
one can assert the unboundless of the multiplicity function $N(y, f,
V\setminus\{x_0\})$ in a neighborhood $V$ of a point $x_0.$ We will
show that Theorem~\ref{th2}, generally speaking, cannot be improved,
in particular: 1) there are examples of open discrete mappings with
finite length $(p, q)$-dis\-tor\-ti\-on with $p\geqslant n,$ that
have an isolated essential singularity, while their multiplicity
function $N(y, f, D\setminus\{x_0\})$ is bounded; 2) the result of
Theorem~\ref{th2} applies only to points $y_0$ in the cluster set
$C(x_0, f)$ but does not extend to arbitrary points $y_0\in {\Bbb
R}^n.$ In other words, there are examples of open discrete mappings
of finite length $(p, q)$-distortion that have an isolated essential
singularity, for which the function $K_{I, p}\left(y, f^{\,-1},
D\right)$ is locally integrable in some neighborhood of some (and
even arbitrary points) $y_0$ in ${\Bbb R}^n\setminus C(x_0, f)$; 3)
It follows from Theorem~A that, under the conditions of this
theorem, $C(x_0, f)=\overline{{\Bbb R}^n}.$ However, under the
conditions of Theorem~\ref{th2}, a similar conclusion cannot,
generally speaking, be drawn. Below, we construct an example of
mapping to which these remarks apply.
\end{remark}

\medskip
\begin{remark}\label{rem2}
As follows from our recent result in~\cite{DS}, if an open discrete
mapping satisfies relation~(\ref{eq5.4AA}) with $K_{I, p}\left(y,
f^{\,-1}, D\right)$ integrable, then an isolated boundary point
$x_0$ for this mapping is always removable. As a corollary, we
obtain the assertion: {\it if a point $x_0$ is essentially
singularity for such a mapping, then the function $K_{I, p}\left(y,
f^{\,-1}, D\right)$ is not integrable.} However, in
Theorem~\ref{th2} we assert something more: the function $K_{I,
p}\left(y, f^{\,-1}, D\right)$ is not only not integrable, but is
not even locally integrable at any point $y_0$ in the cluster set
$C(x_0, f).$ Therefore, the result we present and prove, contained
in the assertion of Theorem~\ref{th2}, is stronger than \cite{DS} in
this respect. On the other hand, the class of mappings involved
in~\cite{DS} is slightly broader than the class of mappings
satisfying~(\ref{eq5.4AA}). Thus, the result of Theorem~\ref{th2}
cannot be deduced from~\cite{DS}, and the result in~\cite{DS} is
also not a consequence of Theorem~\ref{th2}. These results have some
common part, but neither is a generalization of the other.
\end{remark}

\medskip
\begin{remark}
If $p=n,$ then $K_{I, n}\left(y, f^{\,-1},
D\setminus\{x_0\}\right)=:K_I\left(y, f^{\,-1},
D\setminus\{x_0\}\right)\geqslant 1$ for a.e. $y\in
f(D\setminus\{x_0\}).$ Indeed, $K_{I, n}\left(y, f^{\,-1},
D\setminus\{x_0\}\right):=\sum\limits_{x\in D\setminus\{x_0\}\cap
f^{\,-1}(y)}K_{O, n} (x,f),$ however,  $K_{O, n} (x,f)\geqslant 1$
a.e., see e.g. \cite[item~1.7.6]{Fe}. Thus, the requirement
$Q(y)\geqslant c$ for some $c>0$ and for a.e $y\in
f(D\setminus\{x_0\})$ in Theorem~\ref{th2} may be omitted for $p=n.$
\end{remark}

\section{Preliminaries}

Let us recall some definitions. Everywhere below, the boundary and
the closure of the set $A\subset {\Bbb R}^n$ are denoted by
$\partial A$ and $\overline{A},$ respectively, and they should be
understood in the sense of the extended Euclidean space
$\overline{{\Bbb R}^n}={\Bbb R}^n\cup\{\infty\}.$ Also, everywhere
below $D$ is a domain in~${\Bbb R}^n,$ $n\geqslant 2.$

Let $y_0\in {\Bbb R}^n,$ $0<r_1<r_2<\infty$ and
\begin{equation}\label{eq1**}
A(y_0, r_1,r_2)=\left\{ y\,\in\,{\Bbb R}^n:
r_1<|y-y_0|<r_2\right\}\,.\end{equation}
Given sets $E,$ $F\subset\overline{{\Bbb R}^n}$ and a domain
$D\subset {\Bbb R}^n$ we denote by $\Gamma(E,F,D)$ the family of all
paths $\gamma:[a,b]\rightarrow \overline{{\Bbb R}^n}$ such that
$\gamma(a)\in E,\gamma(b)\in\,F$ and $\gamma(t)\in D$ for $t \in (a,
b).$ If $f:D\rightarrow {\Bbb R}^n$ is a given mapping, $y_0\in
f(D)$ and $0<r_1<r_2<d_0=\sup\limits_{y\in f(D)}|y-y_0|,$ then by
$\Gamma_f(y_0, r_1, r_2)$ we denote the family of all paths $\gamma$
in $D$ such that $f(\gamma)\in \Gamma(S(y_0, r_1), S(y_0, r_2),
A(y_0,r_1,r_2)).$ Let $Q:{\Bbb R}^n\rightarrow [0, \infty]$ is a
Lebesgue measurable function.

A mapping $f:D\rightarrow {\Bbb R}^n$ is called {\it discrete} if
the image $\{f^{-1}\left(y\right)\}$ of any point $y\,\in\,{\Bbb
R}^n$ consists of isolated points, and {\it open} if the image of
any open set $U\subset D$ is an open set in ${\Bbb R}^n.$

\medskip
Later, in the extended space $\overline{{{\Bbb R}}^n}={{\Bbb
R}}^n\cup\{\infty\}$ we use the {\it spherical (chordal) metric}
$h(x,y)=|\pi(x)-\pi(y)|,$ where $\pi$ is a stereographic projection
of $\overline{{{\Bbb R}}^n}$ onto the sphere
$S^n(\frac{1}{2}e_{n+1},\frac{1}{2})$ in ${{\Bbb R}}^{n+1},$ namely,
$$h(x,\infty)=\frac{1}{\sqrt{1+{|x|}^2}}\,,$$
\begin{equation}\label{eq3C}
\ \ h(x,y)=\frac{|x-y|}{\sqrt{1+{|x|}^2} \sqrt{1+{|y|}^2}}\,, \ \
x\ne \infty\ne y
\end{equation}
(see, e.g., \cite[Definition~12.1]{Va}).
Further, for the sets $A, B\subset \overline{{\Bbb R}^n}$ let us put
\begin{equation}\label{eq1A}
h(A, B)=\inf\limits_{x\in A, y\in B}h(x, y)\,,\quad
h(A)=\sup\limits_{x, y\in A}h(x ,y)\,, \end{equation}
where $h$ is the chordal distance defined in ~(\ref{eq3C}).

\medskip Let us first formulate a simple but
very important topological statement, which is repeatedly used later
(see, e.g., \cite[Theorem~1.I.5.46]{Ku}).

\medskip
\begin{proposition}\label{pr2}
{\it\, Let $A$ be a set in a topological space $X.$ If $C$ is
connected and $C\cap A\ne \varnothing\ne C\setminus A,$ then $C\cap
\partial A\ne\varnothing.$}
\end{proposition}

\medskip
Let $D\subset {\Bbb R}^n,$ $f:D\rightarrow {\Bbb R}^n$ be a discrete
open mapping, $\beta: [a,\,b)\rightarrow {\Bbb R}^n$ be a path, and
$x\in\,f^{\,-1}(\beta(a)).$ A path $\alpha: [a,\,c)\rightarrow D$ is
called a {\it maximal $f$-lifting} of $\beta$ starting at $x,$ if
$(1)\quad \alpha(a)=x\,;$ $(2)\quad f\circ\alpha=\beta|_{[a,\,c)};$
$(3)$\quad for $c<c^{\prime}\leqslant b,$ there is no a path
$\alpha^{\prime}: [a,\,c^{\prime})\rightarrow D$ such that
$\alpha=\alpha^{\prime}|_{[a,\,c)}$ and $f\circ
\alpha^{\,\prime}=\beta|_{[a,\,c^{\prime})}.$ If $\beta:[a,
b)\rightarrow\overline{{\Bbb R}^n}$ is a path and if
$C\subset\overline{{\Bbb R}^n},$ we say that $\beta\rightarrow C$ as
$t\rightarrow b,$ if the spherical distance $h(\beta(t),
C)\rightarrow 0$ as $t\rightarrow b$ (see \cite[section~3.11]{MRV}),
where $h(\beta(t), C)$ is defined in~(\ref{eq1A}). The following
assertion holds (see~\cite[Lemma~3.12]{MRV}).

\medskip
\begin{proposition}\label{pr3}
{\it Let $f:D\rightarrow {\Bbb R}^n,$ $n\geqslant 2,$ be an open
discrete mapping, let $x_0\in D,$ and let $\beta: [a,\,b)\rightarrow
{\Bbb R}^n$ be a path such that $\beta(a)=f(x_0)$ and such that
either $\lim\limits_{t\rightarrow b}\beta(t)$ exists, or
$\beta(t)\rightarrow \partial f(D)$ as $t\rightarrow b.$ Then
$\beta$ has a maximal $f$-lifting $\alpha: [a,\,c)\rightarrow D$
starting at $x_0.$ If $\alpha(t)\rightarrow x_1\in D$ as
$t\rightarrow c,$ then $c=b$ and $f(x_1)=\lim\limits_{t\rightarrow
b}\beta(t).$ Otherwise $\alpha(t)\rightarrow \partial D$ as
$t\rightarrow c.$}
\end{proposition}

\medskip
Let us proceed to the implementation of the first stage. The
following statement holds, see~\cite[Lemma~2.2]{SevSkv},
cf.~\cite[Theorem~10.12]{Va}.

\medskip
\begin{lemma}\label{lem2}
{\it\, Let $D$ be a domain in ${\Bbb R}^n,$ $n\geqslant 2,$ and
$x_0\in D.$ Then, for any $P>0$ and any neighborhood $U$ of the
point $x_0$ there is a neighborhood $V\subset U$ of the same point,
such that the inequality $M(\Gamma(E, F, D))>P$ holds for any
continua $E, F\subset D$ which intersect $\partial U$ and $\partial
V.$}
\end{lemma}

\medskip
\begin{remark}\label{rem1}
Since the modulus of the family of paths passing through a fixed
point is zero (see section~7.9 in \cite{Va}),  the statement of
Lemma~\ref{lem2} remains true for the case when $x_0$ is an isolated
point of the boundary of the domain $D.$
\end{remark}

\section{On mappings with inverse Poletsky inequality}

We say that {\it $f$ satisfies inverse Poletsky inequality relative
to $p$-modulus} at the point $y_0\in \overline{f(D)},$ $p\geqslant
1,$ if the ratio
\begin{equation}\label{eq2*A}
M_p(\Gamma_f(y_0, r_1, r_2))\leqslant \int\limits_{f(D)\cap
A(y_0,r_1,r_2)} Q(y)\cdot \eta^p (|y-y_0|)\, dm(y)
\end{equation}
holds for any Lebesgue measurable function $\eta:
(r_1,r_2)\rightarrow [0,\infty ]$ such that
\begin{equation}\label{eqA2}
\int\limits_{r_1}^{r_2}\eta(r)\, dr\geqslant 1\,.
\end{equation}

The following statement is true, cf. \cite[Theorem~1.1]{DS}.

\medskip
\begin{theorem}\label{th1}
{\it\, Let $n\geqslant 2,$ $p\geqslant n,$ let $D$ be a domain in
${\Bbb R}^n,$ $x_0\in D,$ and let $f:D\setminus\{x_0\}\rightarrow
{\Bbb R}^n$ be an open discrete mapping that satisfies the
conditions~(\ref{eq2*A})--(\ref{eqA2}) at any point $y_0\in
\overline{D^{\,\prime}}\setminus\{\infty\},$ where
$D^{\,\prime}:=f(D\setminus\{x_0\}).$

\medskip
Assume that, there is $y_0\in  C(x_0, f)$ and a neighborhood $U$ of
$y_0$ such that:

\medskip
1) for any $0<r_1<r_2<r_0:=d(y_0,
\partial U)$ there is a set
$E\subset[r_1, r_2]$ of positive linear Lebesgue measure such that
the function $Q$ is integrable on $S(y_0, r)$ for any $r\in E$
relative to the $(n-1)$-dimensional Hausdorff measure
$\mathcal{H}^{n-1}$ on $S(y_0, r);$

\medskip
2) there exists $y_k\rightarrow y_0$ such that $N(y_k, f,
D\setminus\{x_0\})<\infty$ for any $k\in {\Bbb N}.$

Then $f$ has a continuous extension
$\overline{f}:D\rightarrow\overline{{\Bbb R}^n},$ the continuity of
which should be understood in the sense of the chordal metric $h$
in~(\ref{eq3C}). The extended mapping $\overline{f}$ is open and
discrete in $D.$}
\end{theorem}

\medskip
\begin{proof} Without loss of generalization, we may consider that
$D$ is a bounded domain.  Let us prove this assertion by
contradiction. Let $z_1, z_2\in C(x_0, f)\cap
\partial D^{\,\prime},$ $z_1\ne z_2.$
We may assume that $z_1\ne \infty.$ Then we may also find some
sequences $x_m, x_m^{\,\prime}\in D\setminus \{x_0\},$
$m=1,2,\ldots,$ such that $x_m, x_m^{\,\prime}\rightarrow x_0$ as
$m\rightarrow\infty,$ and $y_m:=f(x_m)\rightarrow z_1$ and
$y^{\,\prime}_m:=f(x_m^{\,\prime})\rightarrow z_2$ as
$m\rightarrow\infty.$ Let $U$ be a neighborhood of $z_1$ such that
$z_2\not \in U.$ Put $0<r_1<r_2<r_0:=d(y_0,
\partial U).$ By the assumption, there is $y\in B(z_1, r_1)$ such that $N(y, f, D\setminus\{x_0\})
<\infty.$ Since $y_m=f(x_m)\rightarrow z_1$ as $m\rightarrow\infty,$
we may consider that $y_m\in B(z_1, r_1)$ for all $m\in {\Bbb N.}$
We join the points $y_m$ and $y$ by a path $\alpha_m:[0,
1]\rightarrow B(z_1, r_1)$ such that $\alpha_m(0)=y_m,$
$\alpha_m(1)=y.$ Let us join the points $y_m$ and $y^{\,\prime}_m$
by a straight line $r_m(t)=y_m+(y^{\,\prime}_m-y_m)t,$
$-\infty<t<\infty.$ Obviously, the ray $\beta_m(t)=r_m(t)|_{t
\geqslant 1}$ does not intersect $U$ and, consequently, does not
intersect $B(z_1, r_2.)$

\medskip
Let $\gamma_m:[0, c_m)\rightarrow D\setminus\{x_0\}$ be a maximal
$f$-lifting of $\alpha_m$ starting at $x_m$ and let
$\gamma^{\,\prime}_m:[0, d_m)\rightarrow D\setminus\{x_0\}$ be a
maximal $f$-lifting of $\beta_m$ starting at $x_m.$ Both maximal
liftings exist by Proposition~\ref{pr3}. By the same Proposition,
there are two situations: either $\gamma_m(t)\rightarrow x_1\in
D\{x_0\}$ as $t\rightarrow c_m,$ or $\gamma_m(t)\rightarrow
\partial {D\setminus\{x_0\}}$ as $t\rightarrow c_m.$
Similarly, $\gamma^{\,\prime}_m(t)\rightarrow x_1\in
D\{x_0\}$ as $t\rightarrow d_m,$ or
$\gamma^{\,\prime}_m(t)\rightarrow
\partial (D\setminus\{x_0\})$ as $t\rightarrow d_m.$

\medskip
Let us to prove that, the situation
$\gamma^{\,\prime}_m(t)\rightarrow x_1\in D\{x_0\}$ as $t\rightarrow
d_m$ is impossible. Indeed, by Proposition~\ref{pr3} we would have
that $d_m=1$ and $f(\omega_1)=\lim\limits_{t\rightarrow
1-0}\beta_m(t).$ Then, from one side, $f(\omega_1)\in D^{\,\prime}$
by the openness of the mapping $f,$ and on the other hand,
$\beta_m(t)\rightarrow \infty$ as $t\rightarrow 1-0$ by the
construction.  Thus, the situation
$\gamma^{\,\prime}_m(t)\rightarrow x_1\in D\setminus\{x_0\}$ as
$t\rightarrow d_m$ is impossible, as required.

\medskip
Due to the mentioned above, since $\partial
(D\setminus\{x_0\})=\partial D\cup\{x_0\},$ we have the following
six different situations:

\medskip
1) for any $m\in {\Bbb N},$ $h(\gamma_{m}^1(t), \partial
D)\rightarrow 0$ as $t\rightarrow c_{m}-0$ and $h(\gamma_{m}^2(t),
\partial D)\rightarrow
0$ as $t\rightarrow d_{m}-0;$

\medskip
2) for any $m\in {\Bbb N},$ $h(\gamma_{m}^1(t), \partial
D)\rightarrow 0$ as $t\rightarrow c_{m}-0,$ but there exists $m_0\in
{\Bbb N}$ such that $h(\gamma_{m_0}^2(t), x_0)\rightarrow 0$ as
$t\rightarrow d_{m_0}-0;$

\medskip
3) there exists $k_0\in {\Bbb N}$ such that $h(\gamma_{k_0}^1(t),
x_0)\rightarrow 0$ as $t\rightarrow c_{k_0}-0,$ in addition, for any
$m\in {\Bbb N},$ $h(\gamma_{m}^2(t),
\partial D)\rightarrow 0$ as $t\rightarrow d_{m}-0;$

\medskip
4) there are numbers $k_0, m_0\in {\Bbb N}$ such that
$h(\gamma_{k_0}^1(t), x_0)\rightarrow 0$ as $t\rightarrow c_{k_0}-0$
and $h(\gamma_{m_0}^2(t), x_0)\rightarrow 0$ as $t\rightarrow
d_{m_0}-0.$

\medskip
5) there exists $k_0\in {\Bbb N}$ such that $h(\gamma_{k_0}^1(t),
p_m)\rightarrow 0$ as $t\rightarrow c_{k_0}-0$ for some $p_n\in
D\setminus\{x_0\},$ in addition, $h(\gamma_{m}^2(t),
\partial D)\rightarrow
0$ as $t\rightarrow d_{m}-0;$

\medskip
6) there exists $k_0\in {\Bbb N}$ such that $h(\gamma_{k_0}^1(t),
p_m)\rightarrow 0$ as $t\rightarrow c_{k_0}-0$ for some $p_m\in
D\setminus\{x_0\},$ in addition, there exists $m_0\in {\Bbb N}$ such
that $h(\gamma_{m_0}^2(t), x_0)\rightarrow 0$ as $t\rightarrow
d_{m_0}-0.$

\medskip
To complete the proof, we need to examine each of these situations
separately. We will argue using the methodology of proving our
recent results in~\cite{DS}.

\medskip
{\bf 1. Let us assume that possibility~1) is fulfilled):} for any
$m\in {\Bbb N},$ $h(\gamma_{m}^1(t), \partial D)\rightarrow 0$ as
$t\rightarrow c_{m}-0$ and $h(\gamma_{m}^2(t),
\partial D)\rightarrow
0$ as $t\rightarrow d_{m}-0.$

\medskip
Put $P>0.$ Let $U:=B(x_0, \delta_0/2),$ and let $V$ be a
neighborhood of $x_0$ which corresponds to Lemma~\ref{lem2} and
Remark~\ref{rem1}. Since by the assumption $x_m, x_m^{\,\prime}\in
D\setminus \{x_0\},$ $m=1,2,\ldots,$ we may find a number $m_0\in
{\Bbb N}$ such that $x_m, x^{\,\prime}_{m}\in V$ for any $m\geqslant
m_0.$
Observe that, for $m\geqslant m_0$
\begin{equation}\label{eq10A}
|\gamma^1_m|\cap \partial V\ne\varnothing, \quad |\gamma^2_m|\cap
\partial V\ne\varnothing\,.
\end{equation}
Indeed, $x_m\in |\gamma^1_m|,$ $x^{\,\prime}_{m}\in |\gamma^2_m|,$
and therefore $|\gamma^1_m|\cap V\ne\varnothing\ne |\gamma^2_m|\cap
V$ for $m\geqslant m_0.$ Besides that, ${\rm diam}\,V\leqslant {\rm
diam}\,U=\delta_0$ and, since $d(|\gamma^1_m|)\geqslant \delta_0>0$
and $d(|\gamma^2_m|)\geqslant \delta_0>0$ for any $m\in {\Bbb N}.$
Now, by Proposition~\ref{pr2} we obtain the relations~(\ref{eq10A}).
Similarly, we may prove that
\begin{equation}\label{eq11B}
|\gamma^1_m|\cap \partial U\ne\varnothing, \quad |\gamma^2_m|\cap
\partial U\ne\varnothing\,.
\end{equation}
Then, by Lemma~\ref{lem2} and Remark~\ref{rem1}
\begin{equation}\label{eq11C}
M(\Gamma(|\gamma^1_m|, |\gamma^2_m|, D\setminus \{x_0\}))>P\,,\qquad
m\geqslant m_0.
\end{equation}
If $p>n,$ since $D$ is bounded, for any $\rho\in {\rm
adm}\,\Gamma(|\gamma^1_m|, |\gamma^2_m|, D\setminus\{x_0\})$ by the
H\"{o}lder inequality
\begin{equation}\label{eq1B}
M(\Gamma(|\gamma^1_m|, |\gamma^2_m|, D\setminus\{x_0\}))\leqslant
\int\limits_{D}\rho^n(x)\,dm(x)\leqslant
\left(\int\limits_{D}\rho^p(x)\,dm(x)\right)^{\frac{n}{p}}\cdot
(m(D))^{\frac{p-n}{p}}\,.
\end{equation}
Passing in (\ref{eq1B}) to $\inf$ over all $\rho\in {\rm
adm}\,\Gamma,$ we obtain that
\begin{equation}\label{eq1E} M(\Gamma(|\gamma^1_m|, |\gamma^2_m|,
D\setminus\{x_0\}))\leqslant \left(M_p(\Gamma(|\gamma^1_m|,
|\gamma^2_m|, D\setminus\{x_0\}))\right)^{\frac{n}{p}}\cdot
(m(D))^{\frac{p-n}{p}}\,.
\end{equation}
It follows from~(\ref{eq11C}) and~(\ref{eq1E}) that
\begin{equation}\label{eq1F}
M_p(\Gamma(|\gamma^1_m|, |\gamma^2_m|, D\setminus\{x_0\}))\geqslant
P^{\frac{p}{n}}\cdot (m(D))^{\frac{n-p}{n}}\,.
\end{equation}
Let us show that the relation~(\ref{eq11C}) and~(\ref{eq1F}) for
$p=n$ and $p>n,$ respectively, are impossible (in particular, each
of them contradicts with the definition of mapping~$f$
in~(\ref{eq2*A})--(\ref{eqA2})). Let $\Gamma_*=\Gamma(|\alpha_m|,
|\beta_m|, D^{\,\prime}).$ Note that
\begin{equation}\label{eq3D}
\Gamma_*>\Gamma(S(z_1, r_1), S(z_1, r_2), A(z_1, r_1, r_2))\,.
\end{equation}
Indeed, let $\gamma\in \Gamma_*,$ $\gamma:[a, b]\rightarrow {\Bbb
R}^n.$ Since $\gamma(a)\in |\alpha_m|\subset B(z_1, r_1)$ and
$\gamma(b)\in |\beta_m|\subset \overline{{\Bbb R}^n}\setminus B(z_1,
r_1),$ by Proposition~\ref{pr2} we may find $t_1\in (a, b)$ such
that $\gamma(t_1)\in S(z_1, r_1).$ Without loss of generalization,
we may assume that $|\gamma(t)-z_1|>r_1$ for $t>t_1.$ Next, since
$\gamma(t_1)\in B(z_1, r_2)$ and $\gamma(b)\in |\beta_m|\subset
{\Bbb R}^n\setminus B(z_1, r_2),$ by Proposition~\ref{pr2} there is
$t_2\in (t_1, b)$ such that $\gamma(t_2)\in S(z_1, r_2).$  Without
loss of generalization, we may assume that $|\gamma(t)-z_1|<r_2$
when $t_1<t<t_2.$ Therefore, $\gamma|_{[t_1, t_2]}$ is a subpath of
$\gamma$ which belongs to~$\Gamma(S(z_1, r_1), S(z_1, r_2), A(z_1,
r_1, r_2)).$ Therefore, the relation~(\ref{eq3D}) is proved.

Let us establish now that
\begin{equation}\label{eq5A}
\Gamma(|\gamma^1_m|, |\gamma^2_m|, D\setminus\{x_0\})>\Gamma_f(z_1,
r_1, r_2)\,.
\end{equation}
Indeed, if the path $\gamma:[a, b]\rightarrow D\setminus\{x_0\}$
belongs to $\Gamma(|\gamma^1_m|, |\gamma^2_m|, D\setminus\{x_0\}),$
then $f(\gamma)$ belongs to $D^{\,\prime},$ and $f(\gamma(a))\in
|\alpha_m|$ and $f(\gamma(b))\in |\beta_m|,$ that is, $f(\gamma)\in
\Gamma_*.$ Then according to the above proof and due to the
ratio~(\ref{eq3D}) the path $f(\gamma)$ has a subpath
$f(\gamma)^{\,*}:=f(\gamma)|_{[t_1, t_2]},$ $a\leqslant
t_1<t_2\leqslant b,$ which belongs to the family $\Gamma(S(z_1,
r_2), S(z_1, r_1), A(z_1, r_1, r_2)).$ Then
$\gamma^*:=\gamma|_{[t_1, t_2]}$ is a subpath of $\gamma$ and it
belongs to~$\Gamma_f(z_1, r_1, r_2),$ as required.

\medskip
In turn, by~(\ref{eq5A}) we have the following:
$$M_p(\Gamma(|\gamma^1_m|, |\gamma^2_m|, D\setminus\{x_0\}))\leqslant$$
\begin{equation}\label{eq11A}
\leqslant M_p(\Gamma_f(z_1, r_1, r_2))\leqslant \int\limits_{A}
Q(y)\cdot \eta^p (|y-z_1|)\, dm(y)\,,
\end{equation}
where $A=A(z_1, r_1, r_2)$ is defined in~(\ref{eq1**}), and $\eta$
is arbitrary Lebesgue measurable function that satisfies
ratio~(\ref{eqA2}) for $r_1:=r_1$ and $r_2:=r_2.$
Put $\widetilde{Q}(y)=\max\{Q(y), 1\},$
\begin{equation}\label{eq13AA}
I=\int\limits_{r_1}^{r_2}\frac{dt}{t^{\frac{n-1}{p-1}}
\widetilde{q}_{z_1}^{1/(p-1)}(t)}\,,
\end{equation}
where $\omega_{n-1}$ denotes the area of the unit sphere ${\Bbb
S}^{n-1}$ in ${\Bbb R}^n,$ and $\widetilde{q}_{z_1}(t)$ is defined
by the relation
$$
\widetilde{q}_{z_1}(r)=\frac{1}{\omega_{n-1}r^{n-1}}\int\limits_{S(z_1,
r)}\widetilde{Q}(y)\,d\mathcal{H}^{n-1}(y)\,.$$
By the assumption of Theorem~\ref{th1}, there exists a set $E\subset
[r_1, r_2]$ of positive linear measure such that
$\widetilde{q}_{z_1}(t)$ is finite for all $t\in E.$ Therefore,
$I\ne 0$ in~(\ref{eq13AA}). In this case, the function
$\eta_0(t)=\frac{1}{It^{\frac{n-1}{p-1}}
\widetilde{q}_{z_1}^{1/(p-1)}(t)}$ satisfies the
relation~(\ref{eqA2}) for $r_1:=\varepsilon_0$ and
$r_2:=\varepsilon^*_1.$ Substituting this function in the right-hand
side of the inequality~(\ref{eq11A}) and applying Fubini theorem, we
obtain that
\begin{equation}\label{eq14AA}
M_p(\Gamma(|\alpha_m|, |\beta_m|, D\setminus\{x_0\}))\leqslant
\frac{\omega_{n-1}}{I^{p-1}}<\infty\,.
\end{equation}
The relation~(\ref{eq14AA}) contradicts with~(\ref{eq11C}) for $p=n$
and~(\ref{eq1F}) for $p>n$. The above contradiction completes the
consideration of the case~\textbf{1)}.

\medskip
{\bf 2. Let us consider the situation  2):} for any $m\in {\Bbb N},$
$h(\gamma_{m}^1(t),
\partial D)\rightarrow 0$ as $t\rightarrow c_{m}-0,$ but there
exists $m_0\in {\Bbb N}$ such that $h(\gamma_{m_0}^2(t),
x_0)\rightarrow 0$ as $t\rightarrow d_{m_0}-0.$

Since $\gamma_{m_0}^2(t)\rightarrow x_0$ as $t\rightarrow
d_{m_0}-0,$ there is a sequence $t_k\rightarrow d_{m_0}-0$ such that
$\gamma_{m_0}^2(t_k)\rightarrow x_0$ as $k\rightarrow\infty.$ We put
$u_k:=\gamma_{m_0}^2(t_k)$ and $v_k:=f(\gamma_{m_0}^2(t_k)).$ Let
also
\begin{equation}\label{eq15F}
D^{\,k*}:={\gamma_{m_0}^2}|_{[0, t_k]}\,, \quad
D^k:=f\left({\gamma_{m_0}^2}|_{[0, t_k]}\right)\,,\quad
k=1,2,\ldots\,.
\end{equation}
In addition, let a path $\alpha_m,$ $m=1,2,\ldots, $ is defined in
the same way as in case~1).

By the definitions of $\alpha_m$ and $D^{\,m*},$ there exists
$\delta_0>0$ such that $d(\alpha_m)\geqslant \delta_0>0$ and
$d(|D^{\,m*}|)\geqslant \delta_0>0$ for all $m=1,2,\ldots .$

\medskip
Let us fix $P>0.$ Let $U:=B(x_0, \delta_0/2),$ and let $V$ be a
neighborhood of the same point $x_0,$ which corresponds to
Lemma~\ref{lem2} and Remark~\ref{rem1}. Reasoning in the same way as
in case~1), we obtain that
\begin{equation}\label{eq11N}
M(\Gamma(|\alpha_m|, |D^{m\,*}|, D\setminus \{x_0\}))>P\,,\qquad
m\geqslant m_0
\end{equation}
for $p=n.$ If $p>n,$
\begin{equation}\label{eq1S}
M_p(\Gamma(|\alpha_m|, |D^{m\,*}|, D\setminus\{x_0\}))\geqslant
P^{\frac{p}{n}}\cdot (m(D))^{\frac{n-p}{n}}\,.
\end{equation}
On the other hand, similarly to~(\ref{eq14AA}) we may show that
\begin{equation}\label{eq14****A}
M_p(\Gamma(|\alpha_m|, |D^{m\,*}|, D\setminus\{x_0\}))\leqslant
\frac{\omega_{n-1}}{I^{p-1}}<\infty\,.
\end{equation}
The relation~(\ref{eq14****A}) contradicts~(\ref{eq11N}) for $p=n$
and~(\ref{eq1S}) for $p>n$. The resulting contradiction completes
the consideration of the case~2).

\medskip
{\bf 3. Let us consider the situation  3):} there exists $k_0\in
{\Bbb N}$ such that $h(\gamma_{k_0}^1(t), x_0)\rightarrow 0$ as
$t\rightarrow c_{k_0}-0,$ in addition, for any $m\in {\Bbb N},$
$h(\gamma_{m}^2(t),
\partial D)\rightarrow 0$ as $t\rightarrow c_{m}-0.$

\medskip
Since $\gamma_{k_0}^1(t)\rightarrow x_0$ as $t\rightarrow
c_{k_0}-0,$ there is a sequence $t_k\rightarrow c_{k_0}-0,$
$k\rightarrow\infty,$ such that $\gamma_{k_0}^1(t_k)\rightarrow x_0$
as $k\rightarrow\infty.$ Put $v_k:=\gamma_{k_0}^1(t_k)$ and
$z_1=z_1(k):=f(\gamma_{k_0}^1(t_k)).$ Also let
\begin{equation}\label{eq15D}
E^{\,k*}:={\gamma_{k_0}^1}|_{[0, t_k]}\,, \quad
E^k:=f\left({\gamma_{k_0}^1}|_{[0, t_k]}\right)\,,\quad
k=1,2,\ldots\,.
\end{equation}
By the definitions of $E^{\,m*}$ and $\beta_m$, there exists
$\delta_0>0$ such that $d(|E^{\,m*}|)\geqslant \delta_0>0$ and
$d(|\beta_m|)\geqslant \delta_0>0$ for all $m=1,2,\ldots .$

\medskip
Let us fix $P>0.$ Let $U:=B(x_0, \delta_0/2),$ and let $V$ be a
neighborhood of the same point $x_0,$ which corresponds to
Lemma~\ref{lem2} and Remark~\ref{rem1}. Reasoning in the same way as
in case~1), we have that
\begin{equation}\label{eq11O}
M(\Gamma(|E^{m\,*}|, |\beta_m|, D\setminus \{x_0\}))>P\,,\qquad
m\geqslant m_0\,,
\end{equation}
in the case $p=n.$ If $p>n,$ then
\begin{equation}\label{eq1T}
M_p(\Gamma(|E^{m\,*}|, |\beta_m|, D\setminus \{x_0\}))\geqslant
P^{\frac{p}{n}}\cdot (m(D))^{\frac{n-p}{n}}\,,\qquad m\geqslant
m_0\,.
\end{equation}
On the other hand, similarly to~(\ref{eq14AA}) or (\ref{eq14****A}),
\begin{equation}\label{eq14*****A}
M_p(\Gamma(|E^{m\,*}|, |\beta_m|, D\setminus\{x_0\}))\leqslant
M_p(\Gamma_{f}(z^1, r_1, r_2))\leqslant
\frac{\omega_{n-1}}{I^{p-1}}<\infty\,.
\end{equation}
The relation~(\ref{eq14*****A}) contradicts with~(\ref{eq11O}) for
$p=n$ and~(\ref{eq1T}) for $p>n$. The resulting contradiction
completes the consideration case~3).

\medskip
{\bf 4. Let us consider the situation  4):} there are numbers $k_0,
m_0\in {\Bbb N}$ such that $h(\gamma_{k_0}^1(t), x_0)\rightarrow 0$
as $t\rightarrow c_{k_0}-0$ and $h(\gamma_{m_0}^2(t),
x_0)\rightarrow 0$ as $t\rightarrow d_{m_0}-0.$

\medskip
Since $\gamma_{k_0}^1(t)\rightarrow x_0$ as $t\rightarrow
c_{k_0}-0,$ there is a sequence $t_k\rightarrow c_{k_0}-0,$
$k\rightarrow\infty,$ such that $\gamma_{k_0}^1(t_k)\rightarrow x_0$
as $k\rightarrow\infty.$ Put $v_k:=\gamma_{k_0}^1(t_k)$ and
$z_1=z_1(k):=f(\gamma_{k_0}^1(t_k)).$ Also let~(\ref{eq15D}) holds.
By the definitions of $E^{\,m*},$ there exists $\delta_0>0$ such
that $d(|E^{\,m*}|)\geqslant \delta_0>0.$

\medskip Similarly, since
$\gamma_{m_0}^2(t)\rightarrow x_0$ as $t\rightarrow d_{m_0}-0,$
there is a sequence $t_k\rightarrow d_{m_0}-0$ such that
$\gamma_{m_0}^2(t_k)\rightarrow x_0$ as $k\rightarrow\infty.$ We put
$u_k:=\gamma_{m_0}^2(t_k)$ and $v_k:=f(\gamma_{m_0}^2(t_k)).$
Let~(\ref{eq15D}) holds.

By the definition of $D^{\,m*},$ there exists $\delta_1>0$ such that
$d(\alpha_m)\geqslant \delta_1>0$ and $d(|D^{\,m*}|)\geqslant
\delta_0>0$ for all $m=1,2,\ldots .$

\medskip
Let us fix $P>0.$ Let $U:=B(x_0, \delta_*),$ where
$\delta_*=\min\{\delta_0, \delta_1\},$ and let $V$ be a neighborhood
of the same point $x_0,$ which corresponds to Lemma~\ref{lem2} and
Remark~\ref{rem1}. Reasoning in the same way as in case~1), we have
that
\begin{equation}\label{eq11_1}
M(\Gamma(|D^{\,m*}|, |E^{m\,*}|, D\setminus \{x_0\}))>P\,,\qquad
m\geqslant m_0\,,
\end{equation}
in the case $p=n.$ If $p>n,$ then
\begin{equation}\label{eq11_2}
M_p(\Gamma(|D^{\,m*}|, |E^{m\,*}|, D\setminus \{x_0\}))\geqslant
P^{\frac{p}{n}}\cdot (m(D))^{\frac{n-p}{n}}\,,\qquad m\geqslant
m_0\,.
\end{equation}
On the other hand, similarly to~(\ref{eq14AA}) or (\ref{eq14****A}),
\begin{equation}\label{eq14_1}
M_p(\Gamma(|E^{m\,*}|, |\beta_m|, D\setminus\{x_0\}))\leqslant
M_p(\Gamma_{f}(z^1, r_1, r_2))\leqslant
\frac{\omega_{n-1}}{I^{p-1}}<\infty\,.
\end{equation}
The relation~(\ref{eq14_1}) contradicts with~(\ref{eq11_1}) for
$p=n$ and~(\ref{eq11_2}) for $p>n$. The resulting contradiction
completes the consideration case~4).

\medskip
{\bf 5. Let us consider the situation  5):} there exists $k_0\in
{\Bbb N}$ such that $h(\gamma_{k_0}^1(t), p_m)\rightarrow 0$ as
$t\rightarrow c_{k_0}-0$ for some $p_m\in D\setminus\{x_0\},$ in
addition, $h(\gamma_{m}^2(t),
\partial D)\rightarrow
0$ as $t\rightarrow d_{m}-0.$

\medskip
In this case by Proposition~\ref{pr3} $c_m=1$ and
$f(p_0)=\lim\limits_{t\rightarrow 1}\alpha_m(t)=y.$ Since by the
assumption $N(y, f, D\setminus\{x_0\})<\infty,$ there is $r>0$ such
that $d(p_m, x_0)\geqslant r>0$ for any $m\in {\Bbb N}.$ Thus,
$d(|\gamma^1_m|)\geqslant \delta_0>0$ for some $\delta_0>0.$ Now,
repeating the reasonings from the item~1), we obtain that, on the
one side, (\ref{eq11C})--(\ref{eq1F}) hold, and on the other hand,
the relation~(\ref{eq11A}) is true. The obtained contradiction
completes the consideration of this item.

\medskip
{\bf 6. Let us consider the situation  6):} there exists $k_0\in
{\Bbb N}$ such that $h(\gamma_{k_0}^1(t), p_m)\rightarrow 0$ as
$t\rightarrow c_{k_0}-0$ for some $p_m\in D\setminus\{x_0\},$ in
addition, there exists $m_0\in {\Bbb N}$ such that
$h(\gamma_{m_0}^2(t), x_0)\rightarrow 0$ as $t\rightarrow
d_{m_0}-0.$

\medskip
As above, $d(p_m, x_0)\geqslant r>0$ for any $m\in {\Bbb N}.$ Thus,
$d(|\gamma^1_m|)\geqslant \delta_0>0$ for some $\delta_0>0.$
Similarly, since $\gamma_{m_0}^2(t)\rightarrow x_0$ as $t\rightarrow
d_{m_0}-0,$ there is a sequence $t_k\rightarrow d_{m_0}-0$ such that
$\gamma_{m_0}^2(t_k)\rightarrow x_0$ as $k\rightarrow\infty.$ We put
$u_k:=\gamma_{m_0}^2(t_k)$ and $v_k:=f(\gamma_{m_0}^2(t_k)).$
Let~(\ref{eq15D}) holds.

By the definition of $D^{\,m*},$ there exists $\delta_1>0$ such that
$d(|D^{\,m*}|)\geqslant \delta_1>0$ for all $m=1,2,\ldots .$ Now,
similarly to~2), we have (\ref{eq11N})--(\ref{eq1S}) on the one
side, and (\ref{eq14****A}) on the other hand. These relations
contradict each other, and this completes the possibility of a
continuous extension $\overline{f}:D\rightarrow \overline{{\Bbb
R}^n}.$ The openness and discreteness of $\overline{f}$ follows by
\cite[Proposition~2.4]{DS}.~$\Box$
\end{proof}

\medskip
\begin{corollary}\label{cor1}
{\,\it The conclusion of Theorem~\ref{th1} remains true if in this
theorem condition 1), which is involved in this theorem, is replaced
by a simpler and more easily verified condition: $Q\in L^1({\Bbb
R}^n).$ }
\end{corollary}

\medskip
\begin{proof}
By the Fubini theorem (see, e.g., \cite[Theorem~8.1.III]{Sa}) we
obtain that
$$\int\limits_{r_1<|x-x_0|<r_2}Q(x)\,dm(x)=\int\limits_{r_1}^{r_2}
\int\limits_{S(x_0, r)}Q(x)\,d\mathcal{H}^{n-1}(x)dr<\infty\,.$$
This means the fulfillment of the condition of the integrability of
the function $Q$ on the spheres with respect to any subset $E_1$ in
$[r_1, r_2].$
\end{proof}

\section{On removable singularities of $FLD$ mappings}

The following statement is true.

\medskip
\begin{theorem}\label{th3}
{\it\, Let $n\geqslant 2,$ $p\geqslant n,$ let $D$ be a domain in
${\Bbb R}^n,$ $x_0\in D,$ and let $f:D\setminus\{x_0\}\rightarrow
{\Bbb R}^n$ be an open discrete mapping with finite length $(p,
q)$-dis\-tor\-ti\-on of $D$ onto $\overline{D^{\,\prime}}.$ Assume
that, there is $y_0\in  C(x_0, f)\setminus\{\infty\}$ and a
neighborhood $U$ of $y_0$ such that, for any $0<r_1<r_2<r_0:=d(y_0,
\partial U)$ there is a set
$E\subset[r_1, r_2]$ of positive linear Lebesgue measure such that
the function $Q:=K_{I, p}\left(y, f^{\,-1}, E\right)$ for $y\in
D^{\,\prime},$ $Q(y)\equiv 0$ for $y\not\in D^{\,\prime},$ is
integrable on $S(y_0, r)$ for any $r\in E$ relative to the
$(n-1)$-dimensional Hausdorff measure $\mathcal{H}^{n-1}$ on $S(y_0,
r).$

Let $Q(y)\geqslant c$ for some $c>0$ and for a.e $y\in {\Bbb R}^n.$
Then $f$ has a continuous extension
$\overline{f}:D\rightarrow\overline{{\Bbb R}^n},$ the continuity of
which should be understood in the sense of the chordal metric $h$
in~(\ref{eq3C}). The extended mapping $\overline{f}$ is open and
discrete in $D.$}
\end{theorem}

\medskip
\begin{proof}
By Remark~\ref{rem2} $f$ satisfies~(\ref{eq2*A}) with $Q(y)=K_{I,
p}\left(y, f^{\,-1}, D\right)$ for $y\in f(D),$ $Q(y)\equiv 0$ for
$y\not\in f(D).$ Assume the contrary, i.e., $x_0$ is an essential
singularity for $f.$ Now, there is  by Theorem~\ref{th1}, $N(y, f,
D\setminus\{x_0\})=\infty$ for any $y\in U$ and some a neighborhood
$U$ of $y_0.$ Observe that,
$$Q:=K_{I, p}\left(y, f^{\,-1}, D\right)=\sum\limits_{x\in E\cap
f^{\,-1}(y)}K_{O, p} (x,f)\geqslant \sum\limits_{x\in E\cap
f^{\,-1}(y)}K_{O, p} (x,f)\geqslant c \sum\limits_{x\in E\cap
f^{\,-1}(y)}1=\infty\,.$$ Thus, the requirement of integrability of
the function $Q$ by the spheres $S(y_0, r)$ for any $r\in E$ does
not hold. The contradiction obtained above proves the
theorem.~$\Box$

\end{proof}

\medskip
\begin{corollary}\label{cor2}
{\,\it The conclusion of Theorem~\ref{th3} remains true if in this
theorem the condition of the integrability of $Q$ under spheres is
replaced by the condition: $Q\in L^1({\Bbb R}^n).$ }
\end{corollary}

\medskip
{\it The proof of Corollary~\ref{cor2} is similar to the proof of
Corollary~\ref{cor1}.}

\medskip
{\bf Proof of Theorem~\ref{th2}} directly follows from
Theorem~\ref{th3}.

\section{An example}

\begin{example} Following~\cite[Example~5.3]{SSD},
cf.~\cite[Proposition~6.3]{MRSY$_2$}, we put $p\geqslant 1$ such
that $n/p(n-1)<1.$ Put $\alpha\in (0, n/p(n-1)).$ Let $f: {\Bbb
B}^n\setminus\{0\} \rightarrow B(0, 2),$
$$f(x)\,=\, \frac{1+|x|^{\,\alpha}}{|x|}\cdot x\,.
$$
Now, $f^{\,-1}(y)=\frac{(|y|-1)^{1/\alpha}}{|y|}\cdot y,$
$f^{\,-1}:B(0, 2)\rightarrow{\Bbb B}^n\setminus\{0\}.$ Obviously,
$f$ is a homeomorphism such that $f\in W_{\rm loc}^{1, n}$ and
$f^{\,-1}\in W_{\rm loc}^{1, n}.$ Thus, $f$ is of finite length $(n,
n)$-distortion, or simply, finite length distortion (see e.g.
\cite[Theorem~8.1]{MRSY$_2$}). Observe that, $C(0, f)={\Bbb
S}^{n-1},$ therefore $f$ has no a continuous extension at the
origin. Using the approach applied
in~\cite[Proposition~6.3]{MRSY$_2$}, we may directly calculate that
$Q:=K_I(f^{\,-1}, y, {\Bbb
B}^n\setminus\{0\})=\frac{|y|}{\alpha(|y|-1)}.$

\medskip
It is immediately clear that the function $Q$ mentioned above is non
integrable in $B(0, 2).$ In addition, Theorem~\ref{th3} implies that
$Q$ cannot be integrable on almost all spheres centered at arbitrary
point $y_0\in {\Bbb S}^{n-1}.$

\medskip
On the other hand, we see that the function $Q$ is integrable in the
neighborhood of any point $y_0\in {\Bbb R}^n,$ if this point does
not belong to the cluster set $C(0, f)={\Bbb S}^{n-1}.$ In
particular, this function is locally integrable over spheres
centered at the points $y_0\in {\Bbb R}^n\setminus C(0, f)={\Bbb
R}^n\setminus {\Bbb S}^{n-1}.$ Therefore, the assertion of
Theorem~\ref{th1}, generally speaking, does not extend to points
$y_0\in {\Bbb R}^n\setminus C(x_0, f).$ We also see that the
multiplicity $N(y, f, {\Bbb B}^n\setminus\{0\})$ is identically
equal to one, since $f$ and $f^{\,-1}$ are homeomorphisms.
Therefore, the conclusion of Theorem~\ref{th1} that the function
$N(y, f, D\setminus\{x_0\})$ is unbounded is, generally speaking,
incorrect unless the additional condition $K_{I, p}\left(y,
f^{\,-1}, D\right)\leqslant Q(y)\cdot N(y, f, D\setminus\{x_0\})$
specified in this theorem is assumed.

\medskip
Finally, $C(0, f)\ne \overline{{\Bbb R}^n},$ so the classical
analogue of the Sokhotski-Casorati-Weierstrass theorem is not true
for this mapping. To summarize, if $Q(y):=K_{I, p}\left(y, f^{\,-1},
D\setminus\{x_0\}\right)$ is simply integrable (or integrable on
almost all concentric spheres), then the weakened version of this
theorem, Theorem~\ref{th1}, cannot be improved in general. That is,
the result of this theorem cannot, in general, be extended to points
that do not belong to the corresponding cluster set, and this set
itself may not be all or ``almost all'' Euclidean space.
Consideration of other classes of functions $Q$ in the context of
this problem will be the subject of our further research.
\end{example}

\medskip
{\bf \noindent Victoria Desyatka} \\
Zhytomyr Ivan Franko State University,  \\
40 Velyka Berdychivska Str., 10 008  Zhytomyr, UKRAINE \\
victoriazehrer@gmail.com

\medskip
\medskip
{\bf \noindent Evgeny Sevost'yanov} \\
{\bf 1.} Zhytomyr Ivan Franko State University,  \\
40 Velyka Berdychivska Str., 10 008  Zhytomyr, UKRAINE \\
{\bf 2.} Institute of Applied Mathematics and Mechanics\\
of NAS of Ukraine, \\
19 Henerala Batyuka Str., 84 116 Slov'yansk,  UKRAINE\\
esevostyanov2009@gmail.com

\end{document}